\documentclass[12pt]{article}
\usepackage{amsxtra,amssymb,amsthm,amsmath,latexsym}

\theoremstyle{plain}
\newtheorem{theorem}{Theorem}

\newtheorem{lemma}{Lemma}
\newtheorem{remark}{Remark}

\newcommand{\refT}[1]{Theorem~\ref{T:#1}}
\newcommand{\refS}[1]{Section~\ref{S:#1}}

\newcommand{\refL}[1]{Lemma~\ref{L:#1}}

\def\bysame{\rule{.5in}{.005in},\ }

\def\lra{\longrightarrow}

\def\dota{\dot{a}}

\def\doth{\dot{h}}
\def\dotu{\dot{u}}
\def\dotw{\dot{w}}

\def\oH{{\overset{\circ}{H}}}
\def\oH1{{\overset{\circ}{H}\kern-.02in{}^1}}

\def\ve{{\varepsilon}}
\def\bee{\begin{equation*}}
\def\eee{\end{equation*}}
\def\be{\begin{equation}}
\def\ee{\end{equation}}

\begin{document}
\title{ Discrepancy principle for DSM}

\author{A.G. Ramm\\
 Mathematics Department, Kansas State University, \\
 Manhattan, KS 66506-2602, USA\\
ramm@math.ksu.edu}

\date{}
\maketitle\thispagestyle{empty}

\begin{abstract}

Let $Ay=f$, $A$ is a linear operator in a Hilbert space $H$,
$y\perp N(A):=\{u:Au=0\}$, $R(A):=\{h:h=Au,u\in D(A)\}$
is not closed, $\|f_\delta-f\|\leq\delta$.
Given $f_\delta$, one wants to construct $u_\delta$ such that
$\lim_{\delta\to 0}\|u_\delta-y\|=0$.
A version of the DSM (dynamical systems method) for finding $u_\delta$
consists of solving the problem 
\bee
  \dotu_\delta(t)=-u_\delta(t)+T^{-1}_{a(t)} A^\ast f_\delta, 
\quad  u(0)=u_0, \eqno{(\ast)}\eee
where $T:=A^\ast A$, $T_a:=T+aI$, and $a=a(t)>0$, $a(t)\searrow 0$
as $t\to\infty$ is suitably chosen.
It is proved that $u_\delta:=u_\delta(t_\delta)$ has the property
$\lim_{\delta\to 0}\|u_\delta-y\|=0$.
Here the stopping time $t_\delta$ is defined by the discrepancy
principle:
\bee
  \int^t_0 e^{-(t-s)}a(s)\|Q^{-1}_{a(s)}f_\delta\| ds=c\delta,
  \eqno{(\ast\ast)}\eee
$c\in(1,2)$ is a constant. Equation $(\ast)$ defines $t_\delta$ uniquely
and $\lim_{\delta\to 0}t_\delta=\infty$.

Another version of the discrepancy principle is proved in this paper:
let $a_\delta$ be the solution to the equation
$a\|Q^{-1}_a f_\delta\|=c\delta$, and $t_\delta$ be the solution to the
equation $a(t)=a_\delta$. If $\lim_{t\to\infty} \frac{\dota(t)}{a^2(t)}=0$,
then $\lim_{\delta\to 0} \|u_\delta-y\|=0$, where 
$u_\delta:=u_\delta(t_\delta)$ and $u_\delta(t)$ solves $(\ast)$.

\end{abstract}

\section{Introduction.}\label{S:1}
Let $A$ be a linear bounded operator in a Hilbert space $H$ 
(or in a Banach space $X$), and equation
 \be\label{e1}  Au=f \ee
be solvable, possibly non-uniquely. Let $N(A)=N$ and $R(A)$
denote the null-space and the range of $A$, respectively. 
Denote by $y$ the (unique) minimal-norm solution to \eqref{e1},
$y\perp N$. Given $f_\delta$, $\|f_\delta-f\|\leq\delta$,
one wants to find a stable approximation $u_\delta$ to $y$:
 \be\label{e2} \lim_{\delta\to 0}\|u_\delta-y\|=0. \ee
There are many ways to do this: variational regularization, 
quasisolutions, iterative regularization 
(see, e.g., \cite{I}, \cite{R1}, \cite{T}).

Here we study a version of the dynamical systems method (DSM) 
for finding $u_\delta$:
 \be\label{e3}
 \dotu_\delta(t)=-u_\delta(t)+T^{-1}_{a(t)} A^\ast f_\delta, 
 \qquad u_\delta(0)=u_0,  \ee
where $T:=A^\ast A$ is selfadjoint, $T_a:=T+aI$, $I$ is the identity 
operator,
  \be\label{e4}
   0<a(t);\quad a(t)\searrow 0\hbox{\ as\ }t\to\infty;
   \quad \lim_{t\to\infty}\frac{\dota}{a}=0,
   \quad \dota:=\frac{da}{dt}.  \ee

The element $u_\delta$ in \eqref{e2} is $u_\delta(t_\delta)$, where 
$u_\delta(t)$ is the solution to \eqref{e3}, and $t_\delta$, the
stopping time, is found from the following equation for the 
unknown $t$:
 \be\label{e5}
 \int^t_0 e^{-(t-s)} a(s) \|Q^{-1}_{a(s)} f_\delta\|ds = c\delta,
 \qquad c\in(1,2), \ee
where $Q:=AA^\ast$ is selfadjoint, $Q_a=Q+aI$, $c$ is a constant, and 
$\|f_\delta\|>c\delta$. This equation we call a discrepancy principle.
About other versions of discrepancy principles see \cite{M},
\cite{R1}-\cite{R5}.

The main result of this paper is the following theorem.

\begin{theorem}\label{T:1}
Assume that \eqref{e4} holds and 
  \be\label{e6}
  \lim_{t\to\infty} e^t a(t) \|Q^{-1}_{a(t)} f_\delta\|=\infty.\ee
Then equation \eqref{e5} has a unique solution $t_\delta$,
  \be\label{e7}  \lim_{\delta\to 0} t_\delta=\infty, \ee
and \eqref{e2} holds with $u_\delta:=u_\delta(t_\delta)$.
\end{theorem}

\begin{remark}\label{R:1}
Assumption \eqref{e6} is always satisfied if $f_\delta\not\in R(A)$.
Indeed, 
  \bee
  \lim_{a\to 0}\|aQ^{-1}_a f_\delta\|^2
  =\lim_{a\to 0} \int^{\|Q\|}_0 
    \frac{a^2d(E_\lambda f_\delta,f_\delta)}{(a+\lambda)^2}
  =\|P f_E\|^2>0, \eee
where $E_\lambda$ is the resolution of the identity, corresponding to the 
selfadjoint operator $Q$, $P$ is the orthoprojector onto the null space 
$N(Q)$ of $Q$, $N(Q)=N(AA^\ast)=N(A^\ast):=N^\ast$,
and $\|P_{N^\ast} f_\delta\| >0$ if $f_\delta\not\in R(A)$,
because $\overline{R(A)}=(N^\ast)^\perp$.
\end{remark}

In \refS{2} we prove \refT{1} and \refT{2}, which says that \eqref{e2} 
holds without assumption \eqref{e6} but with an extra assumption 
$\lim_{t\to\infty}\frac{\dot{a}(t)}{a^2(t)}=0$.

\section{Proofs.}\label{S:2}
Let $h(t):=a(t) \|A^{-1}_{a(t)} f_\delta\|:=a(t)g(t)$.

\begin{lemma}\label{L:1}
Assume \eqref{e6}. Then
  \be\label{e8}
  \lim_{t\to\infty} \frac{e^t h(t)}{\int^t_0 e^s h(s) ds}=1 \ee
provided that
  \be\label{e9}   \lim_{t\to\infty} \frac{\doth}{h}=0. \ee
\end{lemma}

\begin{proof}[Proof of \refL{1}.]
Apply L'Hospital's rule to \eqref{e8} and conclude that \eqref{e8}
holds if
  \be\label{e10}
  \lim_{t\to\infty} \frac{e^t h(t)+e^t\doth}{e^t h(t)} =1. \ee

Relation \eqref{e10} holds if \eqref{e9} holds.

\refL{1} is proved.
\end{proof}

\begin{lemma}\label{L:2}
Relation \eqref{e9} holds if and only if
  \be\label{e11}
  \lim_{t\to\infty} \frac{(h^2)^\bullet}{h^2}=0, \quad 
(h^2)^\bullet:=\frac {d(h^2)}{dt}. \ee
\end{lemma}

\begin{proof}[Proof of \refL{2}.]
We have 
  \be\label{e12}   \lim_{t\to\infty} \frac{(h^2)^\bullet}{h^2}= 
2\lim_{t\to\infty} \frac{\dot{h}}{h}. \ee

\refL{2} is proved.
\end{proof}

We have $h=ag$. Therefore
  \be\label{e13}
  \frac{(h^2)^\bullet}{h^2}
  =\frac{(a^2)^\bullet g^2+a^2(g^2)^\bullet}{a^2 g^2} 
  = \frac{(a^2)^\bullet}{a^2} + \frac{(g^2)^\bullet}{g^2}. \ee
If \eqref{e4} holds, then Lemma 2 implies
  \be\label{e14}   \lim_{t\to\infty} \frac{(a^2)^\bullet}{a^2}=0. \ee
Let us prove that
  \be\label{e15}  \lim_{t\to\infty} \frac{(g^2)^\bullet}{g^2}=0. \ee
Let $E_s$ be the resolution of the identity of the selfadjoint operator 
$Q$. Then
  \be\label{e16}
  g^2(t) = \int^b_0 \frac{d\rho(s)}{[s+a(t)]^2},
  \quad \rho(s):=(E_s f_\delta,f_\delta), \quad b:=\|Q\|, \ee
and
  \be\label{e17}
  \begin{aligned}
  (g^2)^\bullet &= -2\dota \int^b_0 \frac{d\rho(s)}{[s+a(t)]^3}
  =2\frac{|\dota|}{a} \int^b_0 \frac{a(t)d\rho(s)}{[s+a(t)]^3}\\
  & \quad \leq 2 \frac{|\dota(t)|}{a(t)} \int^b_0 
    \frac{d\rho(s)}{[s+a(t)]^2}
  =\frac{2|\dota(t)|}{a(t)} g^2(t), 
  \end{aligned}\ee
where we have used the monotone decay of $a(t)$, which implies
$-\dota=|\dota|$. Our results remain valid for $b=\infty$ and their proofs
are essentially the same.

From \eqref{e4}, \eqref{e16} and \eqref{e17} it follows that
  \be\label{e18}
  0\leq \lim_{t\to\infty} \frac{(g^2)^\bullet}{g^2}\leq 2\lim_{t\to\infty} 
\frac{|\dota|}{a}=0. \ee

We have proved the following lemma.

\begin{lemma}\label{L:3}
If \eqref{e4}  holds, then
  \be\label{e19}  \lim_{t\to\infty} \frac{(h^2(t))^\bullet}{h^2(t)}=0 \ee
\end{lemma}

From Lemmas 1-3 we obtain the following result.

\begin{lemma}\label{L:4}
If \eqref{e4} and \eqref{e6} hold, then
  \be\label{e20}  \lim_{t\to\infty} 
  \frac{\int^t_0 e^{-(t-s)} h(s)ds}{h(t)}=1. \ee
\end{lemma}

\begin{proof}[Proof of \refL{4}.]
Relation \eqref{e20} is equivalent to
  \be\label{e21} \lim_{t\to\infty}
  \frac{e^t h(t)}{\int^t_0 e^s h(s) ds}=1. \ee
By Lemmas 1-3, relation \eqref{e21} holds if
  \be\label{e22} \lim_{t\to\infty} \frac{\dota}{a}=0. \ee
This relation holds by assumption \eqref{e4}.

\refL{4} is proved.
\end{proof}

Now let us prove \refT{1}.

\begin{proof}[Proof of \refT{1}.]
By \refL{4} equation \eqref{e5} for large $t$ can be written as
  \be\label{e23}
  a(t) \|Q^{-1}_{a(t)} f_\delta\|
  =c\delta[1+o(1)], \qquad t\to\infty. \ee
Denote $c[1+o(1)]$ by $c_1$. If $c\in(1,2)$, then $c_1\in(1,2)$ for 
sufficiently large $t$. Fix a sufficiently large $t$ and denote $a(t):=a$. 
Consider the equation 
  \be\label{e24}
  \psi(a):=a^2\|Q^{-1}_a f_\delta\|^2
  =\int^b_0 \frac{a^2 d\rho(s)}{(a+s)^2} =c^2_1 \delta^2. \ee
The function $\psi(a)$ is defined on $(0,\infty)$, is monotonically 
growing and continuous, 
$$\psi(+\infty) = \int^b_0 d\rho= \|f_\delta\|^2>c^2_1\delta^2, \quad
\psi(+0)=\|P f_\delta\|^2,$$ 
where $P$, $P:=E_0-E_{0-},$  is the orthogonal projection onto the 
null-space $N(Q)$ of the operator $Q$, 
$N(Q)=N(AA^\ast)=N(A^\ast):=N^\ast$.
Since $\|f_\delta-f\|\leq\delta$, and $f\in R(A)\perp N(A^\ast)$, 
we have 
  \be\label{e25}
  \|Pf_\delta\|\leq \|P(f_\delta-f)\| +\|Pf\|\leq\delta, \ee
because $\|P\|\leq 1$ and $\|Pf\|=0$. Therefore equation \eqref{e24} 
has a unique solution $a=a_\delta$, 
$\lim_{\delta\to 0}a_\delta=0$. Consequently, equation \eqref{e5}
has a unique solution $t_\delta$, which can be found from the equation
  \be\label{e26}    a(t)=a_\delta. \ee
Since $a(t)\searrow 0$, we have
  \be\label{27}   \lim_{\delta\to 0}t_\delta=\infty. \ee
Finally, let us prove that the element 
$u_\delta:=u_\delta(t_\delta)$ satisfies 
relation \eqref{e2}. First we give a simple proof of 
this relation under an 
additional assumption \eqref{e39}, see below. Then we give a more 
complicated proof which does not use any additional assumptions.

We have
  \be\label{e28}
  u_\delta(t)=e^{-t} u_0+\int^t_0 e^{-(t-s)} T^{-1}_{a(s)}
  A^\ast f_\delta ds. \ee
Thus
  \be\label{e29}
  \begin{aligned}
  \|u_\delta(t)-y\|\leq e^{-t} \|u_0\|
  &+\| \int^t_0 e^{-(t-s)} T^{-1}_{a(s)} A^\ast(f_\delta-f)ds\| \\
  &+\|\int^t_0 e^{-(t-s)} T^{-1}_{a(s)} A^\ast f ds-y\|:= j_1+j_2+j_3.
  \end{aligned} \ee
Clearly, $\lim_{t\to\infty} j_1=0$. We have $\|T^{-1}_a A^\ast\|
\leq\frac{1}{2\sqrt{a}}$, $a>0$.

Indeed, $T^{-1} A^\ast=A^\ast Q^{-1}_a=U Q^{\frac{1}{2}} Q^{-1}_a$, 
where $U$ is an isometry, and we have used the formula
$T^{-1}_a A^\ast=A^\ast Q^{-1}_a$ and the polar decomposition
$A^\ast=UQ^{\frac{1}{2}}$. 
Thus 
$\|T^{-1}A^\ast\|\leq\|Q^{\frac{1}{2}} Q^{-1}_a\|
=\sup_{\delta>0} \frac{s^{\frac{1}{2}}}{s+a} =\frac{1}{2\sqrt{a}}$.
Consequently,
  \be\label{e30}  j_2\leq\frac{\delta}{2\sqrt{a(t)}}.   \ee
Furthermore, $f=Ay$, so
  \be\label{e31}
  \|\int^t_0 e^{-(t-s)} T^{-1}_{a(s)} A^\ast Ay ds-y\|\leq \|y\| e^{-t} 
  +\int^t_0 e^{-(t-s)} a(s) \|T^{-1}_{a(s)} y\|ds,\ee
and
  \be\label{e32}
  \lim_{a\to 0} a^2\|T^{-1}_a y\|^2
  =\lim_{a\to 0} \int^b_0 \frac{a^2d(F_sy,y)}{(a+s)^2}
  =\|\mathcal{P} y\|^2=0,\ee
where $F_s$ is the resolution of the identity corresponding to 
the selfadjoint operator $T=A^\ast A$, $\mathcal{P}$ is the orthogonal 
projection onto
$N(T)=N(A)$, and $\mathcal{P} y=0$ because $y\perp N$.

It follows from \eqref{e31} and \eqref{e32} that
  \be\label{e33}  \lim_{t\to\infty} j_3=0. \ee
Indeed, if
  \be\label{e34}
  \varphi(s):=a(s) \|T^{-1}_{a(s)}y\|\to 0\hbox{\ as\ }s\to\infty,\ee
then, as follows from Lemma 6 (see below), we have
  \be\label{e35}
  \lim_{t\to\infty} \int^t_0 e^{-(t-s)} \varphi(s)ds=0. \ee
If
  \be\label{e36}
  \lim_{\delta\to 0} \frac{\delta}{\sqrt{a(t_\delta)}}=0, \ee
  \be\label{37}  \lim_{\delta\to 0}t_\delta=\infty,  \ee
then \eqref{e29}, \eqref{e30} and \eqref{e33} imply relation \eqref{e2}
with $u_\delta=u_\delta(t_\delta)$.

Equations \eqref{e23} and \eqref{e26} show that
  \be\label{e38}
  a\|Q^{-1}_af_\delta\|=c_1\delta,
  \qquad 1<c_1<2,\quad a:=a_\delta. \ee

We do not have a proof of the relation \eqref{e36}, and it is not
known, in general, if \eqref{e36} holds without additional assumptions. We 
can prove \eqref{e36}
under the following additional assumption:
  \be\label{e39}  y=T^\gamma v, \qquad 0<\gamma<\frac{1}{2}. \ee

\begin{lemma}\label{L:5}
If \eqref{e39} holds and $a=a_\delta$ solves \eqref{e38}, then 
\eqref{e36} holds.
\end{lemma}

\begin{proof}[Proof of \refL{5}]
We have $f=Ay=AT^\gamma v$, and \eqref{e38} implies
  \be\label{e40}
  \begin{aligned}
  \frac{c_1\delta}{\sqrt{a_\delta}} 
   & \leq \sqrt{a_\delta} \left( \|Q^{-1}_{a_\delta} (f-f_\delta)\|
    +\|Q^{-1}_{a_\delta} AT^\gamma v\|\right) \\
   & \leq \sqrt{a_\delta}
    \left(\frac{\delta}{a_\delta}+\|AT^{-1}_{a_\delta}T^\gamma v\| \right)
   =\frac{\delta}{\sqrt{a_\delta}} +\sqrt{a_\delta} \|v\| \ 
   \|T^{-1}_{a_\delta} T^{\frac{1}{2}+\gamma}\|,
  \end{aligned} \ee
where we have used the obvious estimate $\|Q^{-1}_{a}\|\leq \frac 1 a$,
which holds for any $a>0$, and the commutation relation 
$Q^{-1}_{a}A=AT^{-1}_{a}$.
Therefore
  \be\label{e41}
  (c_1-1)\frac{\delta}{\sqrt{a_\delta}}
   \leq \sqrt{a_\delta}\|v\| \max_{s>0} 
\frac{s^{\frac{1}{2}+\gamma}}{s+a_\delta}
  =c\|v\| a^\gamma_\delta\to 0\hbox{\ as\ }\delta\to 0, \ee
where $c=c(\gamma)>0$ is a constant,  $c_1>1$, 
$\max_{s>0}\frac{s^{\frac{1}{2}+\gamma}}{s+a_\delta}=c(\gamma)a^{\gamma 
-\frac 1 2}$, and $c=c(\gamma)=p^p(1-p)^{1-p}$ with $p:=\gamma +\frac 1 
2$.

\refL{5} is proved.
\end{proof}

However, we want to prove relation \eqref{e2} with 
$u_\delta=u_\delta(t_\delta)$
(where $t_\delta$ solves \eqref{e26}, and $a_\delta$ solves \eqref{e38})
{\it without any extra assumptions of the type} \eqref{e39}.

To accomplish this, we propose the following argument.

Equation \eqref{e38} can be written as
  \be\label{e42}
  c_1\delta=\|a Q^{-1}_a f_\delta\|=\|(a+Q-Q) Q^{-1}_a f_\delta\|
  =\|QQ^{-1}_a f_\delta-f_\delta\|
  = \|AT^{-1}_a A^\ast f_\delta-f_\delta\|. \ee
Let us denote $w_\delta:=T^{-1}_{a_\delta}A^\ast f_\delta$, 
where $a_\delta$ is the solution to equation \eqref{e42}. 
It is proved in \cite[p.22]{R1} that
  \be\label{e43}  \lim_{\delta\to 0} \|w_\delta-y\|=0.  \ee
By formula \eqref{e28} the relation 
  \be\label{e44}  \lim_{\delta\to 0} u_\delta(t_\delta)=y \ee
holds due to the following lemma.

\begin{lemma}\label{L:6}
If
  \be\label{e45}  \lim_{t\to\infty} q(t)=q,  \ee
then
  \be\label{e46}  \lim_{t\to\infty} \int^t_0 e^{-(t-s)} q(s)ds=q. \ee
\end{lemma}

\begin{proof}[Proof of \refL{6}]
Take $\ve>0$ arbitrarily small. Using \eqref{e45}, find
$t(\ve)$ such that 
$$|q(t)-q|\leq\frac{\ve}{4}\quad \hbox { for }\quad t\geq t(\ve).$$
Then
  \be\label{e47}
  I:=| \int^t_0 e^{-(t-s)} q(s)ds\-q|
   \leq |\int^{t(\ve)}_0 e^{-(t-s)} q(s)ds|
  + |\int^t_{t(\ve)} e^{-(t-s)} q(s)ds-q|
  :=I_1+I_2. \ee

Take $q(\ve)$ sufficiently large, so that
  \be\label{e48} I_1<\frac{\ve}{2} \quad\hbox { for } \quad 
t>\tau(\ve),  \ee
and
  \be\label{e49}
  \begin{aligned}
  I_2
  & \leq \int^t_{t(\ve)} e^{-(t-s)} |q(s)-q|ds
  +| \int^t_{t(\ve)} e^{-(t-s)}dsq-q| \\
  & \leq\frac{\ve}{4}(1-e^{-(t-t(\ve))}+|q|e^{-(t-t(\ve))}
  <\frac{\ve}{2}\quad \hbox { for } \quad t>\tau(\ve).
  \end{aligned}\ee
Then
  \be\label{e50}  I\leq \ve   \quad \hbox { if } \quad\qquad t>\tau(\ve).  
\ee

\refL{6} is proved.
\end{proof}

\refT{1} is proved.
\end{proof}

We can simplify our arguments if the following extra assumption on 
$a(t)$ is made in addition to \eqref{e4}:
  \be\label{e51}  \lim_{t\to\infty}\frac{\dota(t)}{a^2(t)}=0. \ee

\begin{theorem}\label{T:2}
Assume \eqref{e4} and \eqref{e51}. Let $a_\delta$ solve \eqref{e38}
and $t_\delta$ solve \eqref{e26}. Then relation \eqref{e2} holds with
$u_\delta:=u_\delta(t_\delta)$, where $u_\delta(t)$ is given in 
\eqref{e28}.
\end{theorem}

\begin{proof}[Proof of \refT{2}]
We have
  \be\label{e52}
  \|u_\delta(t_\delta)-y\|
  \leq \|u_\delta(t_\delta)-w_\delta(t_\delta)\|
  +\|w_\delta(t_\delta)-y\|, \ee
where $w_\delta(t):=T^{-1}_{a(t)} A^\ast f_\delta$.
It is known (\cite{R1}, p. 22) that
  \be\label{e53} \lim_{\delta\to 0} \|w_\delta(t_\delta)-y\|=0. \ee

Let us prove that
  \bee \lim_{\delta\to 0} \|u_\delta(t_\delta)-w(t_\delta)\|=0. \eee
We have
  \be\label{e54}
  u_\delta(t)=e^{-t} u_0
  +\int^t_0 e^{-(t-s)} w_\delta(s)ds, \ee
and
  \be\label{e55}
  \int^t_0 e^{-(t-s)} w_\delta(s)ds
  =e^{-(t-s)} w_\delta(s) \bigg|^t_0 - \int^t_0 e^{-(t-s)}
  \dotw_\delta(s)ds. \ee
From \eqref{e54} and \eqref{e55} we get
  \be\label{e56}
  \|u_\delta(t)-w_\delta(t)\|  \leq c e^{-t} +\int^t_0 e^{-(t-s)}
  \|\dotw_\delta(s)\|ds. \ee
From \eqref{e54} and \refL{6} it follows that
  \bee
  \lim_{\delta\to 0} \|u_\delta(t_\delta)-w_\delta(t_\delta)\|=0, \eee
provided that
  \be\label{e57} \lim_{t\to\infty} \|\dotw_\delta(t)\|=0. \ee

To prove \eqref{e57} we use assumption \eqref{e51} and the following 
estimate:
  \be\label{e58}
  \|\dotw_\delta(t)\|^2
  = \int^{\|T\|}_0 \frac{\dota^2(t)}{[a(t)+\lambda]^4}
    d(F_\lambda A^\ast f_\delta, A^\ast f_\delta) 
  \leq \frac{\dota^2(t)}{a^4(t)} \|A^\ast f_\delta\|^2 
  \underset{t\to\infty}{\lra} 0, \ee
where $F_\lambda$ is the resolution of the identity corresponding to the 
operator $T$ (cf formula (32)).
\refT{2} is proved.
\end{proof}

\begin{remark}\label{R:2}
Assumption \eqref{e51} holds, for example, if 
$a(t)=\frac{c_0}{(c_1+t)^b}$ for sufficiently large $t$,
$t>T$, $c_0,c_1b>0$ are constants, $b\in(0,1)$, and $T>0$ is an arbitrary 
large fixed number.
On the initial finite interval $t\in[0,T]$ the function
$a(t)>0$ may decay arbitrarily fast.
\end{remark}


\end{document}